\documentclass[11pt]{article}

\usepackage{times,amssymb,euscript,latexsym,epsfig,amsmath}
\usepackage{amssymb,amsmath}
\usepackage{graphicx}
\usepackage{natbib}
\usepackage{color}
\usepackage{epstopdf}
\usepackage{hyperref}

\def\drawline#1#2{\raise 2.5pt\vbox{\hrule width #1pt height #2pt}}
\def\spacce#1{\hskip #1pt}
\def\solid{\drawline{24}{.5}\nobreak}
\def\bdash{\hbox{\drawline{4}{.5}\spacce{2}}}
\def\dashed{\bdash\bdash\bdash\bdash\hskip-2pt\nobreak}
\def\bdot{\hbox{\drawline{1}{.5}\spacce{2}}}
\def\dotted{\hbox{\leaders\bdot\hskip 24pt}\hskip-2pt\nobreak}

\def\chndot{\hbox {\drawline{9.5}{.5}\spacce{2}\drawline{1}{.5}\spacce{2}\drawline{9.5}{.5}}\nobreak}

\def\trian{\raise 1.25pt\hbox{$\scriptscriptstyle\triangle$}\nobreak}
\def\solidtrian{\raise 1.25pt
\hbox to 3bp{
\hfill}\nobreak}

\def\plus{\raise 1.25pt \hbox{$\scriptscriptstyle +$}\nobreak\ }
\def\x{\raise 1.25pt \hbox{$\scriptscriptstyle \times$}\nobreak\ }

\newcommand\Rey{\mbox{\textit{Re}}}  

\begin{document}

\title{Reduced order models for control of fluids using the Eigensystem Realization Algorithm}

\author{Zhanhua Ma, Sunil Ahuja, and Clarence W. Rowley \\
\small{\{zma, sahuja, cwrowley\} at princeton.edu} \\
\small{\em Mechanical and Aerospace Engineering,} \\
\small {\em Princeton University, Princeton, NJ 08544, USA.}
}

\maketitle

\begin{abstract}
As sensors and flow control actuators become smaller, cheaper, and more pervasive, the use of feedback control to manipulate the details of fluid flows becomes increasingly attractive.  One of the challenges is to develop mathematical models that describe the fluid physics relevant to the task at hand, while neglecting irrelevant details of the flow in order to remain computationally tractable. A number of techniques are presently used to develop such reduced-order models, such as proper orthogonal decomposition (POD), and approximate snapshot-based balanced truncation, also known as balanced POD.  Each method has its strengths and weaknesses: for instance, POD models can behave unpredictably and perform poorly, but they can be computed directly from experimental data; approximate balanced truncation often produces vastly superior models to POD, but requires data from adjoint simulations, and thus cannot be applied to experimental data.

In this paper, we show that using the Eigensystem Realization Algorithm (ERA) \citep{JuPa-85}, one can theoretically obtain exactly the same reduced order models as by balanced POD.  Moreover, the models can be obtained directly from experimental data, without the use of adjoint information.  The algorithm can also substantially improve computational efficiency when forming reduced-order models from simulation data.  If adjoint information is available, then balanced POD has some advantages over ERA: for instance, it produces modes that are useful for multiple purposes, and the method has been generalized to unstable systems.  We also present a modified ERA procedure that produces modes without adjoint information, but for this procedure, the resulting models are not balanced, and do not perform as well in examples.  We present a detailed comparison of the methods, and illustrate them on an example of the flow past an inclined flat plate at a low Reynolds number.
\end{abstract}

\section{Introduction}
\label{sec:introduciton}
In the last decade, substantial developments have been made in the area of model-based feedback flow control of fluids: for instance, see the recent reviews by \cite{KimBew-07, CaSoWiRoAl-08, ChoiJeonKim08}.  In many applications, the focus is on how to apply actuation to maintain the flow around an equilibrium state of interest, for instance to delay transition to turbulence, or control separation on a bluff body.  Linear control theory provides efficient tools for the analysis and design of feedback controllers.   However, a significant challenge is that models for flow control problems are often very  high dimensional, e.g., on the order of ${\cal O} (10^{5 \sim9})$,  so large that it becomes computationally infeasible to apply linear control techniques.   To address this issue, model reduction, by which a low-order approximate model is obtained, is therefore widely employed.

Several techniques are available for model reduction, many of which involve projection onto a set of modes.  These may be global eigenmodes of a linearized operator \citep{Aakervik_et_al-jfm-07}, modes determined by proper orthogonal decomposition (POD) of a set of data \citep{HLB-96}, and various variants of POD, such as including shift modes \citep{Noack-03}.  An approach that is used widely for model reduction of linear systems is balanced truncation \citep{Moore-81}, and while this method is computationally intractable for systems with very large state spaces (dimension $\gtrsim 10^5$), recently an algorithm for computing approximate balanced truncation from snapshots of linearized and adjoint simulations has been developed \citep{Rowley-ijbc05} and successfully applied to a variety of high-dimensional flow control problems \citep{IlakRowley-pof08, Ahuja:2009, BaBrHe-09}.     In this method, sometimes called balanced POD, one obtains two sets of modes (primal and adjoint) that are bi-orthogonal, and uses those for projection of the governing equations, just as in standard POD.  Compared to most other methods, including POD, balanced truncation has key advantages, such as {\it a priori\/} error bounds, and guaranteed stability of the reduced-order model (if the original high-order system is stable).   Balanced POD is an approximation of exact balanced truncation that is computationally tractable when the number of states is very large (for instance, up to $10^7$), and typically produces models that are far more accurate than standard POD models.  For instance, for a linearized channel flow investigated in~\cite{IlakRowley-pof08}, even though the first 5 POD modes capture over 99.7\% of the energy in a dataset exhibiting large transient growth, a low-dimensional model obtained by projection onto these modes completely misses the transient growth.  By contrast, a 3-mode balanced POD model captures the transient growth nearly perfectly; to do as well with a standard POD model, 17 modes were required.

The main steps of balanced POD include (a) taking snapshots from impulse responses of the linearized and adjoint systems, (b) computing a singular value decomposition (SVD) of a matrix formed from inner products of these snapshots, (c) constructing primal modes and adjoint modes from the resulting singular vectors, and (d) projecting the high-dimensional dynamics onto these modes.

While effective in many examples, balanced POD also faces challenges, especially for use with experimental data.  The main restriction is that balanced POD requires snapshots of impulse-response data from an adjoint system, and adjoint information is not available for experiments.

To address this issue, here we describe an algorithm widely used for system identification and model reduction, the eigensystem realization algorithm (ERA) \citep{JuPa-85}.  This algorithm has been used for problems in fluid mechanics, primarily as a system-identification technique for flow control~\cite{Cattafesta-97, CaKaCoGi-06}, but also for model reduction \cite{Gaitonde-03, Silva-04}. Our main result, presented in Section~\ref{sec:era}, is that, for linear systems, ERA theoretically produces {\em exactly the same} reduced-order models as balanced POD, with no need of an adjoint system, and at {an order of magnitude lower} computational cost.  This result implies that one can realize approximate balanced truncation even in experiments, and can also improve computational efficiency in simulations.  We note that  ERA and snapshot-based approximate balanced truncation have been applied together in a model reduction procedure in~\cite{DjCaMy-08}.  However, the theoretical equivalence between these two algorithms was not explored in that work.

We present a comparison between balanced POD and ERA, and show that if adjoint information is available, balanced POD also has its own advantages. In particular, balanced POD provides sets of bi-orthogonal primal/adjoint modes for the linear system, and can be directly generalized to unstable systems.  In Section~\ref{sec:ERApseudomodes}, we discuss a modified ERA algorithm that, in the absence of adjoint simulations, uses ``pseudo-adjoint modes'' to compute reduced order models; however, this method does not produce balanced models, and performs worse than balanced POD in examples. In Section~\ref{sec:numerical}, we illustrate these methods using a numerical example of the two-dimensional flow past an inclined plate, at a low Reynolds number.

\section{The eigensystem realization algorithm as snapshot-based approximate balanced truncation}
\label{sec:era}
In this section, we summarize the steps involved in approximate balanced truncation (balanced POD), and the Eigensystem Realization Algorithm, and show that they are equivalent.

Balanced truncation involves first constructing a a coordinate transformation that ``balances'' a linear input-output system, in the sense that certain measures of controllability and observability (the Gramian matrices) become diagonal and identical~\citep{Moore-81}.  A reduced-order model is then obtained by truncating the least controllable and observable states, which correspond to the smallest diagonal entries in the transformed system.   Unfortunately, the exact balanced truncation algorithm is not tractable for the large state dimensions encountered in fluid mechanics.  However, an approximate, snapshot-based balanced truncation algorithm, referred to as Balanced Proper Orthogonal Decomposition (balanced POD) was proposed in~\cite{Rowley-ijbc05}, and has been used successfully in several examples~\citep{IlakRowley-pof08, Ahuja:2009, BaBrHe-09}.

The second technique, the eigensystem realization algorithm (ERA), has been used both for system identification and for model reduction, and it is well known that the models produced by ERA are approximately balanced \citep{Ga-96, JuPh-01}. Here we show further that, theoretically,  ERA produces exactly the same reduced order models as balanced POD.   This equivalence indicates that ERA can be regarded as an approximate balanced truncation method, in the sense that, before truncation, it implicitly realizes a coordinate transformation under which a pair of approximate controllability and observability Gramians are exactly balanced.  This feature distinguishes ERA from other model reduction methods that first realize truncations and then balance the reduced order models.  Note that  in ERA the Gramians, and the balancing transformation itself, are never explicitly calculated, as we will also show in the following discussions. 

For both techniques, we will consider a high-dimensional, stable, discrete-time linear system, described by
\begin{equation}
\begin{aligned}
x(k+1) &= Ax(k)+Bu(k)\\
 y(k)&=Cx(k),
\end{aligned}
\label{LTI}
\end{equation}
where $k\in\mathbb{Z}$ is the time step index, $u(k)\in\mathbb{R}^p$ denotes a vector of inputs (for instance, actuators or disturbances), $y(k)\in\mathbb{R}^q$ a vector of outputs (for instance, sensor measurements, or simply quantities that one wishes to model), and $x(k)\in\mathbb{R}^n$ denotes the state variable (for instance, flow variables at all gridpoints of a simulation).  These equations may arise, for instance, by discretizing the Navier-Stokes equations in time and space, and linearizing about a steady solution, as will be demonstrated in the example in Section~\ref{sec:numerical}.  The goal is to obtain an approximate model that captures the same relationship between inputs~$u$ and outputs~$y$, but with a much smaller state dimension:
\begin{equation}
\begin{aligned}
x_{r}(k+1) &= A_{r}x_{r}(k)+B_{r}u(k)\\
 y(k)&=C_{r}x_{r}(k)
\end{aligned}
\label{LTIReduced}
\end{equation}
where the reduced state variable~$x_r(k)\in\mathbb{R}^r$, $r\ll n$.  We consider the discrete-time setting, because we are primarily interested in discrete-time data from simulations or experiments.

\subsection{Snapshot-based approximate balanced truncation (balanced POD)}
\label{subsec:bpod}
Here, we give only a brief overview of the balanced POD algorithm, and for details of the method, we refer the reader to~\cite{Rowley-ijbc05}.  The algorithm involves three main steps:
\begin{itemize}
\item {\bf Step 1:} {\em Collect snapshots.}  Run impulse-response simulations of the primal system (\ref{LTI}) and collect $m_{c}+1$ snapshots of {\em states}~$x(k)$ in $m_{c}P+1$~steps:
 \begin{equation}
X = \begin{bmatrix} B & A^{P}B & A^{2P}B & \cdots & A^{m_{c}P}B \end{bmatrix},
\label{X}
\end{equation}
where $P$ is the sampling period.  In addition, run impulse-response simulations for the adjoint system
\begin{equation}
\begin{aligned}
z(k+1) &= A^{*}z(k)+C^{*}v(k)
\end{aligned}
\label{LTIAdjoint}
\end{equation}
where the asterisk $^{*}$ stands for adjoint of a matrix, and collect $m_{o}+1$ snapshots of states~$z(k)$ in $m_{o}P+1$~steps:
\begin{equation}
Y = \begin{bmatrix} C^{*}& \left(A^{*}\right)^{P}C^{*} & \left(A^{*}\right)^{2P}C^{*} \cdots & \left(A^{*}\right)^{m_{o}P}C^{*} \end{bmatrix}.
\label{Y}
\end{equation}
Calculate the generalized Hankel matrix,
\begin{equation}
H = Y ^{*}X.
\label{HankelXY}
\end{equation}
\item {\bf Step 2:} {\em Compute modes.}  Compute the singular value decomposition of $H$:
\begin{equation}
 H = U\Sigma V^{*}=\begin{bmatrix}U_{1} & U_{2}\end{bmatrix}\begin{bmatrix}\Sigma_{1} & 0 \\ 0 & 0\end{bmatrix}
 \begin{bmatrix}V_{1}^{*}\\V_{2}^{*}\end{bmatrix}
 =U_{1}\Sigma_{1} V_{1}^{*}
\label{SVDH}
\end{equation}
where the diagonal matrix $\Sigma_{1}\in\mathbb{R}^{n_{1}\times n_{1}}$ is invertible and includes all non-zero singular values of $H$, $n_{1}=\text{rank}(H)$, and $U_{1}^{*}U_{1} = V_{1}^{*}V_{1} = I_{n_{1}\times n_{1}}$.  Choose $r\leq n_{1}$. Let $U_{r}$ and $V_{r}$ denote the sub-matrices of $U_{1}$ and $V_{1}$  that include their first  $r$ columns, and $\Sigma_{r}$ the first $r\times r$ diagonal block of $\Sigma_{1}$.  Calculate
\begin{equation}
\begin{aligned}
\Phi_{r} = XV_{r}\Sigma_{r}^{-\frac{1}{2}}; \quad \Psi_{r} = YU_{r}\Sigma_{r}^{-\frac{1}{2}}.
\end{aligned}
\label{modes}
\end{equation}
where the columns of $\Phi_{r}$ and $\Psi_{r}$ are respectively the first $r$ primal and adjoint modes  of system~(\ref{LTI}). The two sets of modes are bi-orthogonal: $\Psi_{r}^{*}\Phi_{r} = I_{r\times r}$.

\item {\bf Step 3:} {\em Project dynamics.}    
The system matrices in the reduced order model~(\ref{LTIReduced}) are
\begin{equation}
\begin{aligned}
A_{r} =\Psi_{r}^{*}A\Phi_{r}; \quad B_{r} =  \Psi_{r}^{*}B; \quad C_{r} =C\Phi_{r}.
\end{aligned}
\label{ReducedABCBPOD}
\end{equation}
\end{itemize}

Note that the $n\times n$ controllability/observability Gramians are approximated by the matrices $XX^*$ and $YY^*$.  The reduced-order model~(\ref{LTIReduced}) is obtained by considering a subspace $x= \Phi_{r}x_{r}$, and projecting the dynamics~(\ref{LTI}) onto this subspace using the adjoint modes given by~$\Psi_{r}$. It was shown in \cite{Rowley-ijbc05} that $\Phi_{r}$ and $\Psi_{r}$ respectively form the first $r$ columns of the balancing transformation/inverse transformation that  \emph{exactly} balance the approximate controllability/observability  Gramians $XX^*$ and~$YY^*$; see more discussion in Section~\ref{sec:ERApseudomodes}.

\subsection{The eigensystem realization algorithm}
\label{subsec:era}
The eigensystem realization algorithm (ERA) was proposed in~\cite{JuPa-85} as a system identification and model reduction technique for linear systems.     
The algorithm follows three main steps \citep{JuPa-85, JuPh-01}:
\begin{itemize}
\item {\bf Step 1:} Run impulse-response simulations/experiments of the system (\ref{LTI})  for $(m_{c} + m_{o})P+2$ steps, where $m_{c}$ and $m_{o}$ respectively reflect how much effect is taken for considering controllability and observability, and $P$ again is the sampling period.  Collect the snapshots of the {\em outputs}~$y$ in the following pattern:
\begin{equation}
\begin{aligned}
&\left(  CB, \quad CAB, \quad  CA^P B, \quad CA^{P+1} B,   \quad \ldots \right. \\
&\left. CA^{m_{c}P}B, \quad CA^{m_{c}P+1}B,  \quad \ldots \quad CA^{(m_{c}+m_{o})P}B, \quad C A^{(m_{c + m_{o}})P +1 } B\right).
\end{aligned}
\label{snapshots}
\end{equation}
The terms $CA^kB$ are commonly called {\em Markov parameters}.  Construct a generalized Hankel matrix $H\in \mathbb{R}^{q(m_{o}+1)\times p(m_{c}+1)}$
\begin{equation}
H = \begin{bmatrix}
CB & CA^{P}B &  \cdots & CA^{m_{c}P}B \\
CA^{P}B & CA^{2P}B &\cdots &  CA^{(m_{c}+1)P}B\\
\vdots&\vdots &\ddots  &\vdots & \\
CA^{m_{o}P}B &  CA^{(m_{o}+1)P}B & \cdots & CA^{(m_{c}+m_{o})P}B
\end{bmatrix}.
\label{Hankel}
\end{equation}

\item {\bf Step 2:} Compute SVD of $H$, exactly as in~(\ref{SVDH}), to obtain $U_{1}$, $V_{1}$, $\Sigma_{1}$. Let $r\leq\text{rank}(H)$.  Let $U_{r}$ and $V_{r}$ denote the sub-matrices of $U_{1}$ and $V_{1}$  that include their first  $r$ columns, and $\Sigma_{r}$ the first $r\times r$ diagonal block of $\Sigma_{1}$. 

\item {\bf Step 3:}  
The reduced $A_{r}$, $B_{r}$ and $C_{r}$ in (\ref{LTIReduced}) are then defined as
\begin{equation}
\begin{aligned}
A_{r} &= \Sigma^{-\frac{1}{2}}_{r}U_{r}^{*}H{'}V_{r}\Sigma_{r}^{-\frac{1}{2}}; \quad\\
B_{r} &= \text{the first $p$ columns of } \Sigma_{r}^{\frac{1}{2}}V_{1}^{*}; \quad\\
C_{r} &= \text{the first $q$ rows of } U_{r}\Sigma_{r}^{\frac{1}{2}}
\end{aligned}
\label{ReducedABC}
\end{equation}
where 
\begin{equation}
H{'}= \begin{bmatrix}
CAB & CA^{P+1}B &\cdots &  CA^{m_{c}P+1}B\\
\vdots&\vdots &\ddots  &\vdots & \\
CA^{m_{o}P+1}B &  CA^{(m_{o}+1)P+1}B & \cdots & CA^{(m_{c}+m_{o})P+1}B
\end{bmatrix},
\label{HankelPrime}
\end{equation}
which can again be constructed directly  from the collected snapshots~(\ref{snapshots}).
\end{itemize}

\subsection{Theoretical equivalence between ERA and balanced POD}
\label{subsec:equivalence}
The first observation is that, with $X$ and $Y$ given by~(\ref{X}) and~(\ref{Y}), the generalized Hankel matrices obtained in balanced POD and ERA, respectively by~(\ref{HankelXY}) and~(\ref{Hankel}), are theoretically identical. The theoretical equivalence between the two algorithms then follows immediately:   First, $H'$ given in (\ref{HankelPrime}) satisfies $H{'}= Y^{*}A X$, which implies the matrices $A_{r}$ obtained in the two algorithms are identical.   To show the equivalence of $B_r$, first note that the  SVD~(\ref{SVDH}) leads to $\Sigma_{1}^{-\frac{1}{2}} U_{1}^{*}H =  \Sigma_{1}^{\frac{1}{2}}V_{1}^{*}$, which, by definition of $U_{r}$, $V_{r}$, $\Sigma_{r}$, implies $\Sigma_{r}^{-\frac{1}{2}} U_{r}^{*}H =  \Sigma_{r}^{\frac{1}{2}}V_{r}^{*}$.
(Note that it does {\emph {not}} imply $H = U_{r}\Sigma_{r}V_{r}^{*}$, since $U_{r}U_{r}^{*}$ is not the identity.)
Thus, in balanced POD,
$B_{r} = \Psi_{r}^{*}B = \Sigma_{r}^{-\frac{1}{2}} U_{r}^{*}Y^{*}B$, which equals the first $p$ columns of $\Sigma_{r}^{-\frac{1}{2}} U_{r}^{*}H = \Sigma_{r}^{\frac{1}{2}}V_{r}^{*}$, which is the $B_{r}$ given by ERA.  Similarly, the  SVD~(\ref{SVDH}) leads to $HV_{1}\Sigma_{1}^{-\frac{1}{2}} =  U_{1}\Sigma_{1}^{\frac{1}{2}}$ and then  $HV_{r}\Sigma_{r}^{-\frac{1}{2}} =  U_{r}\Sigma_{r}^{\frac{1}{2}}$.  Thus, in balanced POD, $C_{r} = C\Phi_{r} = CXV_{r}\Sigma_{r}^{-\frac{1}{2}}$, which equals the first $q$ rows of $HV_{r}\Sigma_{r}^{-\frac{1}{2}} = U_{r}\Sigma_{r}^{\frac{1}{2}}$, the $C_{r}$ given by ERA.    In summary, we have:

\paragraph{\bf Main result.}  {\it The reduced system matrices $A_{r}$, $B_{r}$ and $C_{r}$ generated in balanced POD and ERA, respectively by~(\ref{ReducedABCBPOD}) and~(\ref{ReducedABC}), are theoretically identical}.

In practice, these two algorithms may generate slightly different reduced order models, because the Hankel matrices calculated in the two algorithms are usually not exactly the same, due to small numerical inaccuracies in adjoint simulations, and/or in matrix multiplications needed to compute the sub-blocks in the Hankel matrices.  In the following discussions, we compare these two algorithms in more detail.

\subsection{Comparison between ERA and balanced POD}
\label{subsec:comparison}
While ERA and balanced POD produce theoretically identical reduced-order models, the techniques differ in several important ways, both conceptually and computationally.  Neither ERA nor balanced POD calculate Gramians explicitly, but balanced POD does construct approximate controllability and observability matrices~$X$ and~$Y^{*}$, from which one calculates the generalized Hankel matrix~$H$ and balancing transformation. Balanced POD thus incurs additional computational cost, because one needs to construct the adjoint system (\ref{LTIAdjoint}), run adjoint simulations for~$Y$, and then calculate each block of~$H$ by matrix multiplication.  Thus we see that the advantages of ERA include:
\begin{enumerate}
\item {\bf Adjoint-free:} ERA is a feasible balanced truncation method for experiments, since it needs only the output measurements from the response to an impulsive input.   Note that ERA has been successfully applied in several flow control experiments \citep{Cattafesta-97, CaKaCoGi-06}, as  a system-identification technique rather than a balanced-truncation method.  In practice, input-output sensor responses are often collected by applying a broadband signal to the inputs, and the {ARMARKOV} method \citep{AkeBer-97b, LimPhaLon-98} can then be used to identify the Markov parameters, or even directly the generalized Hankel matrix, from the input-output data history.

\item {\bf Computational efficiency:} For large problems, typically the most computationally expensive component of computing balanced POD is constructing the generalized Hankel matrix~$H$ in~(\ref{HankelXY}), as this involves computing inner products of all of the (large) primal and adjoint snapshots with each other. ERA is significantly more efficient at constructing the matrix $H$ in~(\ref{Hankel}), since only the first row and last column of block matrices, i.e., the $(m_{c} + m_{o} +1)$ Markov parameters, need be obtained by matrix multiplication.  All the other $m_{c} \times m_{o}$ block matrices in~$H$ are copies of other blocks, and need not be recomputed.  For balanced POD, the matrix~$H$ is obtained by computing all the $(m_{c}+1) \times (m_{o}+1)$ matrix multiplications (inner products) between corresponding blocks in  $Y^{*}$ and $X$ in~(\ref{HankelXY}).  Thus, for example, if $m_{c} = m_{0} = 200$, the computing time  needed for constructing $H$ in ERA will be about only $1\%$ of that in balanced POD.  See Table~\ref{tab:era_savings} for a detailed comparison on computational efficiency between balanced POD and ERA in the example of the flow past an inclined flat plate.   
\end{enumerate}

At the same time, balanced POD also provides its own advantages:
\begin{enumerate}
\item {\bf Sets of bi-orthogonal primal/adjoint modes:} Balanced POD provides sets  of bi-orthogonal primal/adjoint  modes, the columns of $\Phi_{r}$ and $\Psi_{r}$. In comparison, without the adjoint system, ERA cannot provide the primal and adjoint modes.   At best, the primal modes may be computed, using the first equation in~(\ref{modes}), if the matrix $X$~(\ref{X}) is stored (in addition to the Markov parameters).   But the adjoint modes cannot be computed without solutions of the adjoint system.  In this sense, balanced POD incorporates more of the physics of the system (the two sets of bi-orthogonal modes), while ERA is purely based on input-output data of the system.  The primal/adjoint modes together can be useful for system analysis and controller/observer design purposes in several ways:  for instance, in flow control applications,  a large-amplitude region from the most observable mode (the leading adjoint mode) can be a good location for actuator placement.  Also, although balanced POD is a linear method, a nonlinear system can be projected onto these sets of modes to obtain a nonlinear low-dimensional model.  For instance, the transformation $x = \Phi_{r}x_{r}$, $x_{r} = \Psi_{r}^{*}x$ can be employed to reduce a full-dimensional nonlinear model $\dot{x} = f(x)$ to a low-dimensional system $\dot{x}_{r} = \Psi_{r}^{*}f( \Phi_{r}x_{r})$. Finally, if parameters (such as Reynolds number or Mach number) are present in the original equations, balanced POD can retain these parameters in the reduced-order models. When the values of parameters change, the reduced order model by balanced POD may still be valid and perform well; see~\cite{IlakRowley-pof08} for an application to linearized channel flow. 

\item {\bf Unstable systems:} Balanced POD has been extended to neutrally stable \citep{ Ma:2008} and unstable systems \citep{Ahuja:2009}. In those cases, one first calculates the right/left eigenvectors corresponding to the neutral/unstable eigenvalues of the state-transition matrix $A$, using direct/adjoint simulations. Using these eigenvectors, the system is projected onto a stable subspace and then balanced truncation is realized for the stable subsystem.  ERA for unstable systems is still an open problem, if adjoint operators are not available.  However, we note that, once the stable subsystem is obtained, ERA can still be applied to it and efficiently realize its approximate balanced truncation.

\end{enumerate}

\paragraph{ERA for systems with high-dimensional outputs. }
The method of output projection proposed in \cite{Rowley-ijbc05} makes it computationally feasible to realize approximate balanced truncation for systems with high-dimensional outputs---for instance, if one wishes to model the entire state~$x$, say the flow field in the entire computational or experimental domain.  This method involves projecting the outputs onto a small number of POD modes, determined from snapshots of~$y$ from the impulse-response dataset. This method can be directly incorporated into ERA as follows: First, run impulse response simulations of the original system and collect Markov parameters as usual. Then, compute the leading POD modes of the dataset of Markov parameters and stack them as columns of a matrix~$\Theta$. Left multiply those Markov parameters by $\Theta^{*}$ to project the outputs onto these POD modes.  A generalized Hankel matrix is then constructed using these modified Markov parameters, and the usual steps of ERA follow.    

\section{A modified ERA method using pseudo-adjoint modes}
\label{sec:ERApseudomodes}
We have seen that one of the drawbacks of ERA is that it does not provide modes that could be used, for instance, for projection of nonlinear dynamics, or to retain parameters in the models.  More precisely, using ERA, one may still obtain primal modes $\Phi_1= XV_{1}\Sigma_{1}^{-1/2}$ as in balanced POD (see~(\ref{SVDH}--\ref{modes})), as long as the state snapshots are collected and stored in $X$.  But it is not possible to obtain the corresponding adjoint modes $\Psi_1$ necessary for projection, without performing adjoint simulations to gather snapshots for the matrix~$Y$.  This is a severe drawback, as adjoint solutions can be expensive to perform, and are not available for experimental data.  One idea, proposed in~\cite{OrSpeCar-08}, is to define a set of approximate adjoint modes using the Moore-Penrose pseudo-inverse of $\Phi_1$:
\begin{equation}
	\tilde\Psi_1 = \Phi_1 (\Phi_1^*\Phi_1)^{-1}.
	\label{eq:pseudomodes}
\end{equation}
We will call the adjoint modes as defined above the {\em pseudo-adjoint modes} corresponding to the modes~$\Phi_1$.
The system matrices of a $r$-dimensional reduced-order model ($r\leq \text{rank}(H)$) generated by this approach then read
\begin{equation}
\begin{aligned}
A_{r} =\tilde{\Psi}_{r}^{*}A\Phi_{r}; \quad B_{r} =  \tilde{\Psi}_{r}^{*}B; \quad C_{r} =C\Phi_{r},
\end{aligned}
\label{ReducedABCpseudo}
\end{equation}
where $\Phi_{r}$ and $\tilde{\Psi}_{r}$ are respectively the first $n\times r$ sub-blocks of $\Phi_{1}, \tilde{\Psi}_{1}$. 

While this idea does produce a set of modes that can be used for projection, we show in this section that, unfortunately, the resulting transformation is {\em not\/} a balancing transformation, and does not produce models that are an approximation to balanced truncation.  In fact, the resulting models are closer to those produced by the the standard POD/Galerkin method: as with standard POD/Galerkin, the method performs well as long as the most controllable and most observable directions coincide.  However, when these directions differ (as is the case for many problems of interest, including the example in Section~\ref{sec:numerical}), the method performs poorly.  These systems in which controllable and observable directions do not coincide are precisely the systems for which balanced POD and ERA give the most improvement over the more traditional POD/Galerkin approach.

\subsection{Transformed approximate Gramians}
\label{sub:transformed_approximate_gramians}

First, let us recall in what sense the model-reduction procedures described in Section~\ref{sec:era} are approximations to balanced truncation.  Suppose that we have an approximation of the controllability and observability Gramians, factored as
\begin{equation}
	W_c = XX^*,\qquad W_o = YY^*,
	\label{eq:approximate_grams}
\end{equation}
where $X$ and $Y$ are the matrices of snapshots from~(\ref{X}) and~(\ref{Y}).
In balanced POD, we define the primal modes as columns of $\Phi_{1} = XV_{1}\Sigma_{1}^{-\frac{1}{2}}$, and the adjoint modes as columns of  $\Psi_{1} = YU_{1}\Sigma_{1}^{-\frac{1}{2}}$, where $U_1$, $V_1$, and $\Sigma_1$ are defined in~(\ref{SVDH}).  We will assume in this section that the number of columns of $X$ and $Y$ (the number of snapshots, $m_c$ and~$m_o$, respectively) is smaller than the number of rows (the state dimension, $n$), which is always true for the large fluids systems of interest here.

Then balanced POD is an approximation to balanced truncation in the following sense: as shown in the appendix of~\cite{Rowley-ijbc05} (the proof of Proposition 2), we may construct a full (invertible, $n\times n$) transformation
\begin{equation}
	T = \begin{bmatrix}
		\Phi_1 & \Phi_2
	\end{bmatrix}
\end{equation}
by choosing $\Phi_2$ such that $\Psi_1^*\Phi_2=0$. That is, columns of $\Phi_2$ are orthogonal to the adjoint modes, which are columns of~$\Psi_1$.  The inverse transformation then has the form
\begin{equation}
	T^{-1} = \begin{bmatrix}
		\Psi_1^*\\
		\Psi_2^*
	\end{bmatrix}
	\label{eq:tinv_def}
\end{equation}
where $\Psi_1$ is the matrix of adjoint modes, and $\Psi_2$ is defined by~(\ref{eq:tinv_def}).  Then, Proposition~2 of~\cite{Rowley-ijbc05} states that the transformed approximate Gramians (\ref{eq:approximate_grams}) have the form
\begin{equation}
	T^{-1}W_c(T^{-1})^* = \begin{bmatrix}
		\Sigma_1 & 0\\ 0 & M_1
	\end{bmatrix},\qquad
	T^*W_o T = \begin{bmatrix}
		\Sigma_1 & 0\\0 & M_2
	\end{bmatrix},
	\label{eq:xform1}
\end{equation}
and furthermore the product of the approximate Gramians, in the transformed coordinates, is
\begin{equation}
	T^{-1}W_cW_oT = \begin{bmatrix}
		\Sigma_1^2 & 0\\0 & 0
	\end{bmatrix}.
	\label{eq:xform2}
\end{equation}
In this sense, the transformation~$T$ balances the approximate Gramians as closely as possible: the Gramians are block diagonal, and the upper-left blocks are equal and diagonal.  Furthermore, all of the states in the lower-right block (i.e., involving $M_1$ and~$M_2$ above) are either unobservable or uncontrollable, as they do not appear in the product of the Gramians.

However, if the pseudo-adjoint modes $\tilde\Psi_1$ are used in place of the true adjoint modes $\Psi_1$, then this result does not hold, as we now show.  Note that, in order for the first block of rows of $T^{-1}$ to equal $\tilde\Psi_1^*$, we must now define
\begin{equation}
	\tilde T = \begin{bmatrix}
		\Phi_1 & \tilde \Phi_2
	\end{bmatrix}
	\label{eq:Tpseudo_def}
\end{equation}
where $\tilde \Psi_1^{*}\tilde \Phi_2=0$.  Since the range of $\tilde\Psi_1$ equals the range of~$\Phi_1$, this is then equivalent to choosing $\tilde\Phi_2$ such that its columns are orthogonal to the columns of $\Phi_1$ (the {\em primal} modes), while when the ``true'' adjoint modes are used, columns of $\Phi_2$ are chosen to be orthogonal to the {\em adjoint} modes~$\Psi_1$.

Defining $\tilde\Psi_2$ by
\begin{equation}
	\tilde T^{-1} = \begin{bmatrix}
		\tilde\Psi_1^*\\
		\tilde\Psi_2^*
	\end{bmatrix},
\end{equation}
one can then show that, as long as $\operatorname{rank}(X)\le\operatorname{rank}(Y)$\footnote{If $\operatorname{rank}(X)>\operatorname{rank}(Y)$, then the situation is worse, and the transformed controllability Gramian is not block diagonal, nor does its upper-left block equal $\Sigma_1$.}, the transformed Gramians have the form
\begin{equation}
	\tilde T^{-1}W_c(\tilde T^{-1})^* = \begin{bmatrix}
		\Sigma_1 & 0\\ 0 & \tilde M_1
	\end{bmatrix},\qquad
	\tilde T^*W_o \tilde T = \begin{bmatrix}
		\Sigma_1 & M_3\\M_3^* & \tilde M_2
	\end{bmatrix},\qquad
	\tilde T^{-1}W_cW_o\tilde T = \begin{bmatrix}
		\Sigma_1^2 & \Sigma_1 M_3 \\ \tilde M_1 M_3^* & 0
	\end{bmatrix},
\label{eq:xform_pseudo}
\end{equation}
with
\begin{equation}
	M_3 = \Sigma_1 \Psi_1^*\tilde\Phi_2,
\end{equation}
where $\Psi_1=YU_1\Sigma_1^{-1/2}$ are the true adjoint modes.
Note that, when the true adjoint modes are used to define the inverse~(\ref{eq:tinv_def}), then $M_3=0$, since $\Psi_1^*\Phi_2=0$.  However, when pseudo-adjoint modes are used, then $M_3$ is no longer zero, and in fact, can be quite large.

\begin{figure}
\centering
\includegraphics[width=0.8\linewidth]{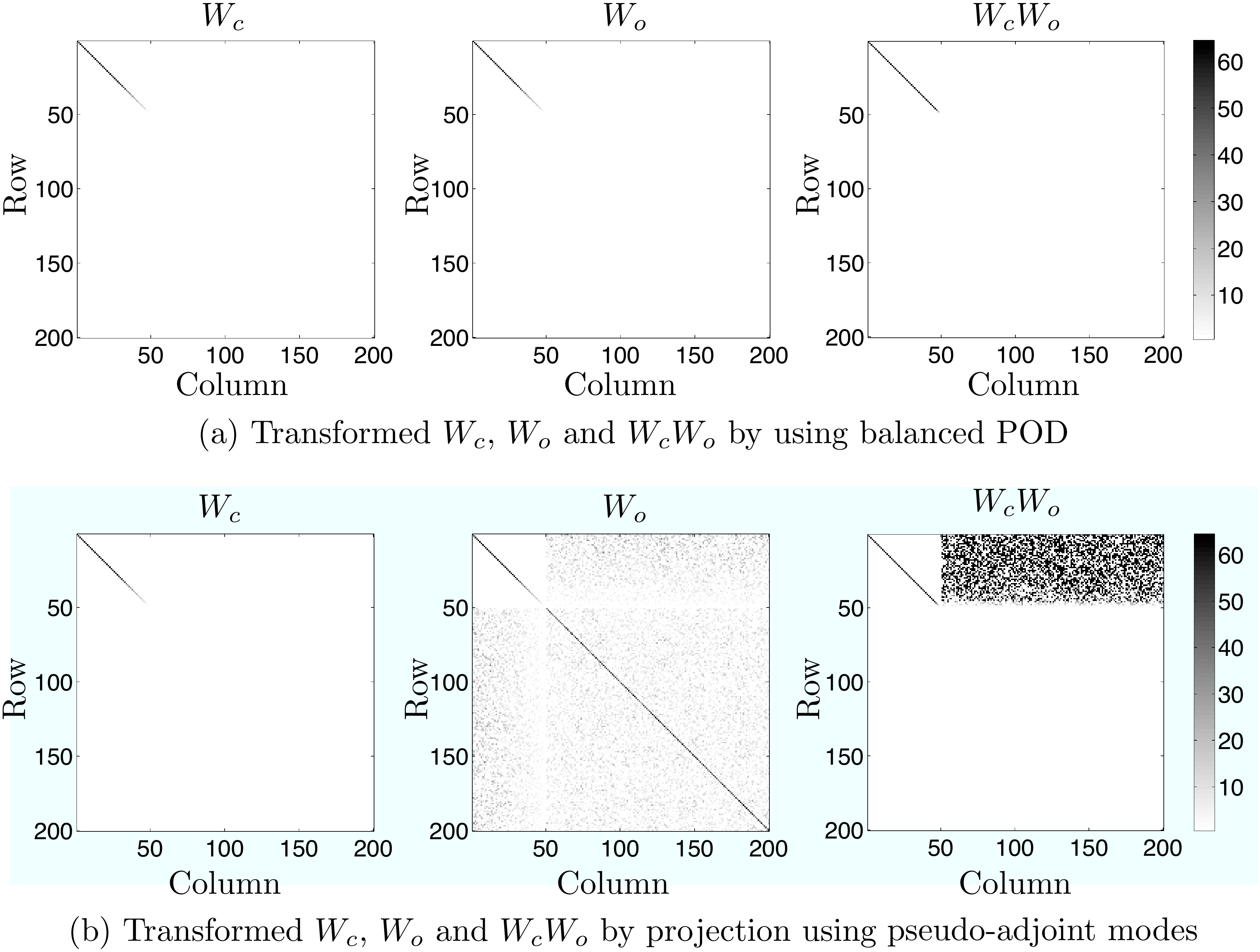}
\caption{Transformed Gramian matrices: (a) using true adjoint modes (eq.~(\ref{eq:xform1}--\ref{eq:xform2})) and (b) using pseudo-adjoint modes (eq.~(\ref{eq:xform_pseudo})).  Here, $X$ and $Y$ in~(\ref{eq:approximate_grams}) are random matrices with $n=200$ states and $m_c=m_o=50$ snapshots.}
\label{fig:random_xform}
\end{figure}

An example is shown in Figure~\ref{fig:random_xform}, which shows the magnitude of the elements of the transformed Gramians, where $X$ and~$Y$ in~(\ref{eq:approximate_grams}) are chosen at random.  Note that when true adjoint modes are used, the transformed Gramians are equal and diagonal, while when the pseudo-adjoint modes are used, the off-diagonal blocks of the transformed observability Gramian, and the product of the Gramians, have significant magnitude.

Thus, when pseudo-adjoint modes are used, the resulting transformation is not, in general, a balancing transformation: even though the upper-left blocks of the transformed Gramians are still equal and diagonal, the transformed observability Gramian is not block diagonal, and so its eigenvalues and eigenvectors do not correspond to those of the transformed controllability Gramian.  Note that this is the whole point of balanced truncation: to transform to coordinates in which the most controllable directions (dominant eigenvectors of~$W_c$) correspond to the most observable directions (dominant eigenvectors of~$W_o$).  Therefore, while the approximate balanced truncation procedure described in Section~\ref{subsec:bpod} exactly balances the approximate Gramians, transforming by pseudo-adjoint modes does not represent balancing in any meaningful sense.

Note that the matrix $M_3$ describes the degree to which projection using pseudo-adjoint modes fails to balance the approximate Gramians.  This matrix equals zero if the adjoint modes (columns of~$\Psi_1$) are spanned by the primal modes (columns of~$\Phi_1$).  However, $M_3$ is the largest when the dominant adjoint modes (columns of $\Psi_1$) are nearly orthogonal to the dominant primal modes (columns of $\Phi_1$).  Unfortunately, this is the case in many problems of interest, in particular those involving non-normality: the directions spanned by the primal modes often do not coincide with the directions spanned by the adjoint modes.

In the next section, we apply this approach to the flow past a flat plate, and compare it to the methods described in Section~\ref{sec:era}.

\section{Example: flow past an inclined flat plate}
\label{sec:numerical}

In this section, we illustrate the application of ERA as an approximate balanced truncation method using a numerical example, by obtaining reduced-order models of a large-dimensional fluid system. We compare the resulting models with those obtained using the balanced POD method of~\cite{Rowley-ijbc05}, ERA with pseudo-adjoint modes as described in Section~\ref{sec:ERApseudomodes}, and the standard POD/Galerkin method~\cite{HLB-96}.

\subsection{Model problem and parameters}
\label{subsec:parameters}

\begin{figure}
\centering
\includegraphics[height=1.3in]{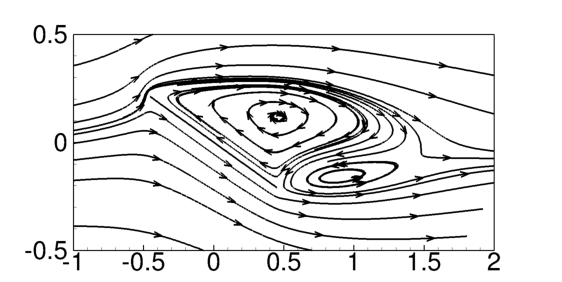}
\includegraphics[height=1.3in]{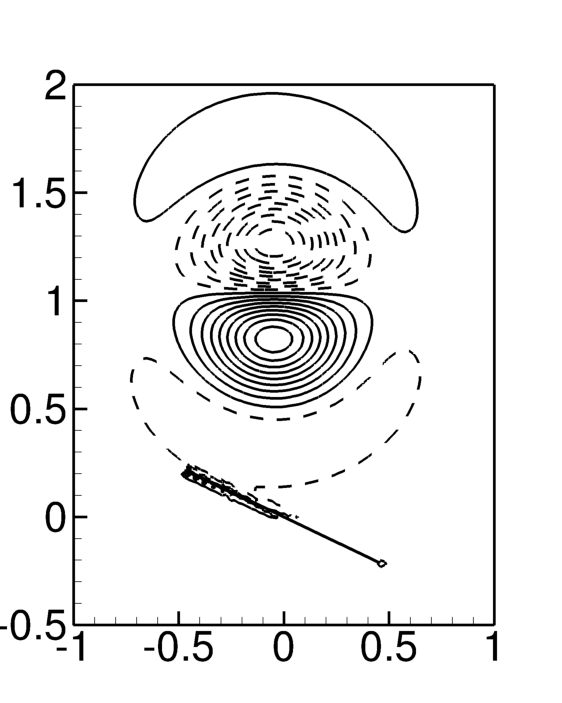}
\caption{Streamlines of the stable steady state past a flat plate at~$\alpha = 25^\circ$~(left), and the contour-lines of the vorticity field obtained from an impulsive input to the actuator~(right). }
\label{figure:flowplate}
\end{figure}

The model problem that we consider is a two-dimensional uniform flow over a flat plate inclined at an angle $\alpha=25^\circ$, at a low Reynolds number~$\Rey = 100$. At these conditions, the flow asymptotically reaches a stable steady state, the streamlines of which are plotted in Figure~\ref{figure:flowplate}. The numerical method used for all computations is a fast formulation of the immersed boundary method developed by~\cite{ColTai-07}, and solves for the vorticity field at each time step. We treat farfield boundary conditions using the multiple-grid scheme described in~\cite{ColTai-07} (Section~4) with five nested grids, each with $250\times 250$ points.  The finest grid covers the region~$[-2,3]\times[-2.5,2.5]$, and the largest grid covers the region~$[-32,48]\times[-40,40]$, where lengths are non-dimensionalized by the chord of the flat plate, whose center is located at the origin. The time step used for all simulations is~$0.01$ (nondimensionalized by chord and freestream velocity).   The numerical model is the same as that considered in~\cite{Ahuja:2009}  where balanced POD is applied for feedback controller design to stabilize an {\em unstable} steady state corresponding to a high angle of attack. However, here we consider the case of a {\em stable} steady state (with an angle of attack  at $25^\circ$), for comparison of reduced order models. 

\subsection{Input and output}
\label{subsec:inputoutput}
The governing equations are first linearized about the stable steady state, resulting in a high-dimensional model in the form of equation~(\ref{LTI}), where the state~$x$ consists of the discrete vorticity field at the grid points. See~\cite{Ahuja:2009} for the details of the linearized (and adjoint) equations and their numerical formulations. The system input $u$ is a disturbance (or actuator) shown in Figure~\ref{figure:flowplate}, modeled as a localized body force in the vicinity of the leading edge. 
We consider the output to be the entire velocity field: this is important for capturing the flow physics, and is often needed to represent cost functions used in optimal control design. Since the output is very high-dimensional, in ERA and balanced POD reduction procedures we use output projection described at the end of  Section~\ref{sec:era}, projecting the velocity field onto the leading POD modes of the velocity snapshots obtained from the impulse response simulation.

\subsection{Reduced-order models}
\label{subsec:rom}
ERA is applied to the full-dimensional linearized system to construct a reduced-order model. 
With a sampling period of 50 time steps, 400 adjacent pairs of Markov parameters, as in~(\ref{snapshots}), are collected from an impulse response simulation. 
Since these parameters are a projection of the velocity fields onto the leading POD modes, for an output projection of order~$m$, the number of inner products required is~$4m \times 10^2$ for construction of each~$H$ and~$H'$ (see Section~\ref{subsec:comparison}). 

For comparison, balanced POD is also used to compute the same reduced-order models. Adjoint simulations are performed with the POD modes as initial conditions to compute the matrix~$Y$ of~(\ref{Y}). The matrices~$X$ and~$Y$ are assembled by stacking 200 snapshots from the linearized and each of the adjoint simulations, and in turn, the generalized Hankel matrix~$H = Y^\ast X$ is computed. For an output projection of order~$m$, the number of inner products required to compute $H$ is~$4m \times 10^4$, which is 50~times more than that to compute $H$ and $H'$ in total for ERA. 

We also compare reduced-order models using standard POD modes, and ERA with pseudo-adjoint modes, as described in Section~\ref{sec:ERApseudomodes}.  The first   100 primal modes are used to compute the pseudo-adjoint modes..

For the given case, a comparison between the computational cost using ERA and using balanced POD  is shown in Table~\ref{tab:era_savings}.  Results verify that ERA substantially improves computational efficiency in forming reduced-order models.
\begin{table} \centering
\begin{tabular}{l c c}
\hline
Steps in computing & \multicolumn{2}{c}{Approximate time (CPU hours)}  \\
 \cline{2-3}
reduced-order models  & balanced POD  & ERA   \\
\hline
1. Linearized impulse response & 2 & 4 \\
2. Computation of POD modes & 2 & 2 \\
3. Adjoint impulse responses  & 30 & - \\
\hskip0.2in(10 in number) & & \\
4. Computation of the Hankel matrix & 7 & 0.2 \\
5. Singular value decomposition & 0.05 & 0.05 \\
6. Computation of modes & 1 & - \\
7. Computation of models & 0.02 & 0.02\\
\hline
\end{tabular}
\caption[Comparison of the computational costs of approximate balanced truncation and ERA.]{Comparison of the computational times required for various steps of the algorithms using balanced POD and ERA.  The times are given for a 10-mode output projected system. The Hankel matrix is constructed using (a) 200 state-snapshots from each linearized and adjoint simulations for balanced POD, and (b) 400 Markov parameters~(outputs) for ERA.}
\label{tab:era_savings}
\end{table}

\begin{figure}
\centering
\includegraphics[width=0.9\linewidth]{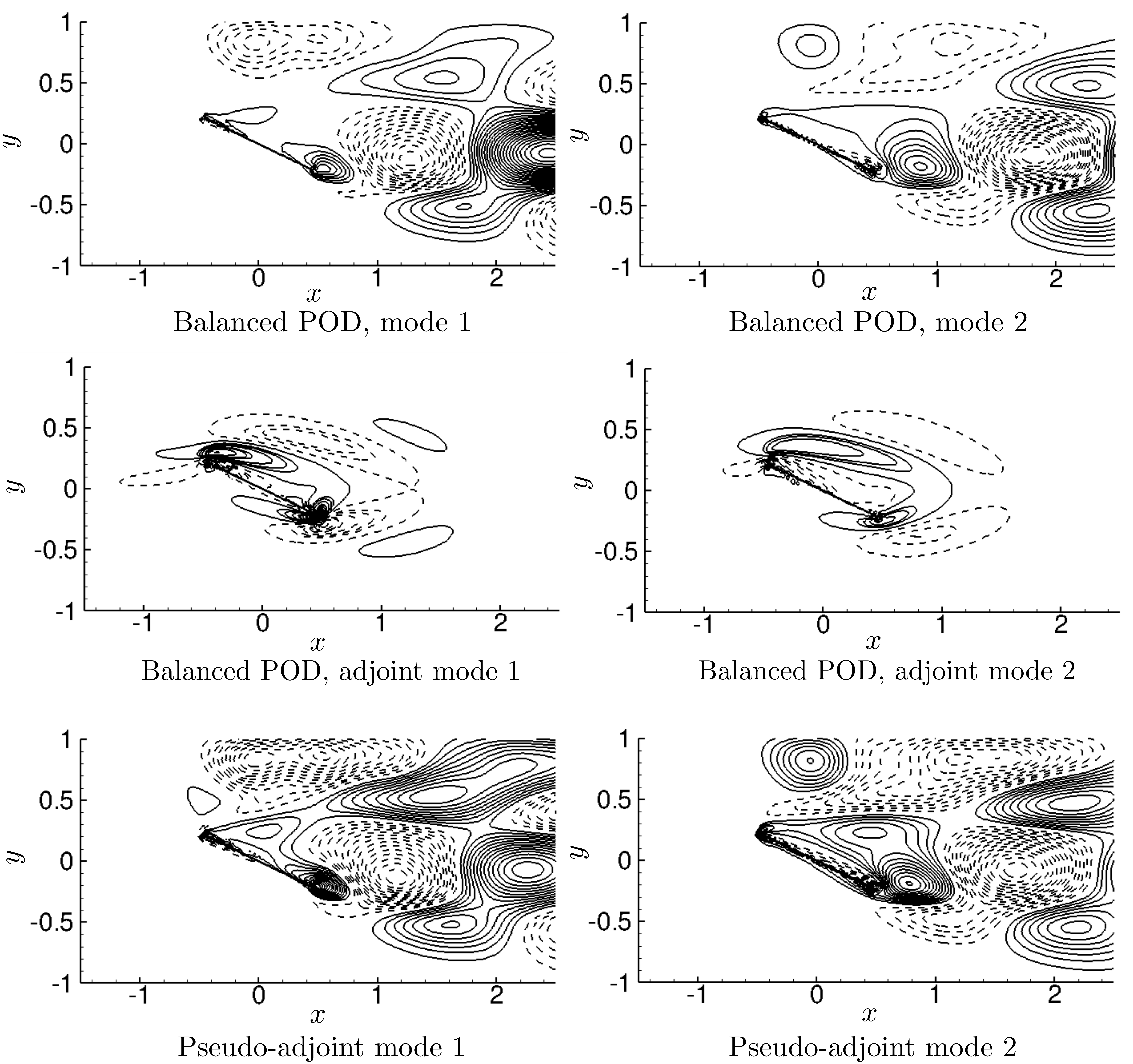}
\caption{ The first two primal and adjoint modes computed by using balanced POD, and the first two pseudo-adjoint modes computed by using~(\ref{eq:pseudomodes}) and the first 100 primal modes. Modes are illustrated using contour plots of the vorticity field.} 
\label{figure:modes}
\end{figure}

Next, we compare the reduced-order models.  Figure~\ref{figure:modes} shows the leading two primal modes and true adjoint modes from balanced POD, compared with the leading two pseudo-adjoint modes. The pseudo-adjoint modes look quite different from the true adjoint modes, and the flow structures actually more closely resemble the leading {\em primal} modes.  This result is not surprising, since the pseudo-adjoint modes are always linear combinations of the snapshots from the primal simulations, while the true adjoint modes are linear combinations of snapshots from adjoint simulations.  Following the discussion in the last section, the poor approximation of the adjoint modes suggests that the pseudo-adjoint modes may produce poor reduced order models for this example, as we will verify below.

\begin{figure}
\centering
\includegraphics[width=2.5in]{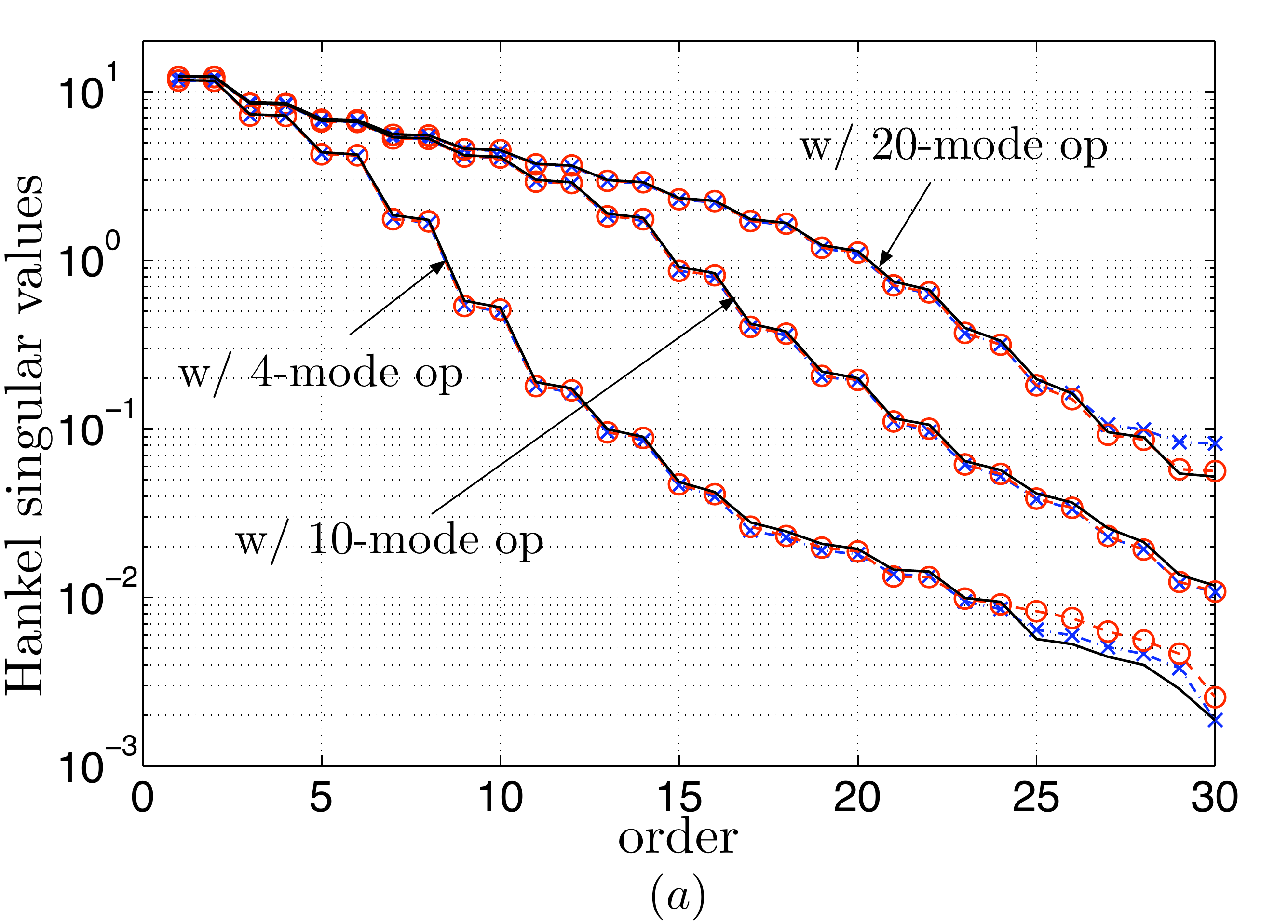}
\includegraphics[width=2.5in]{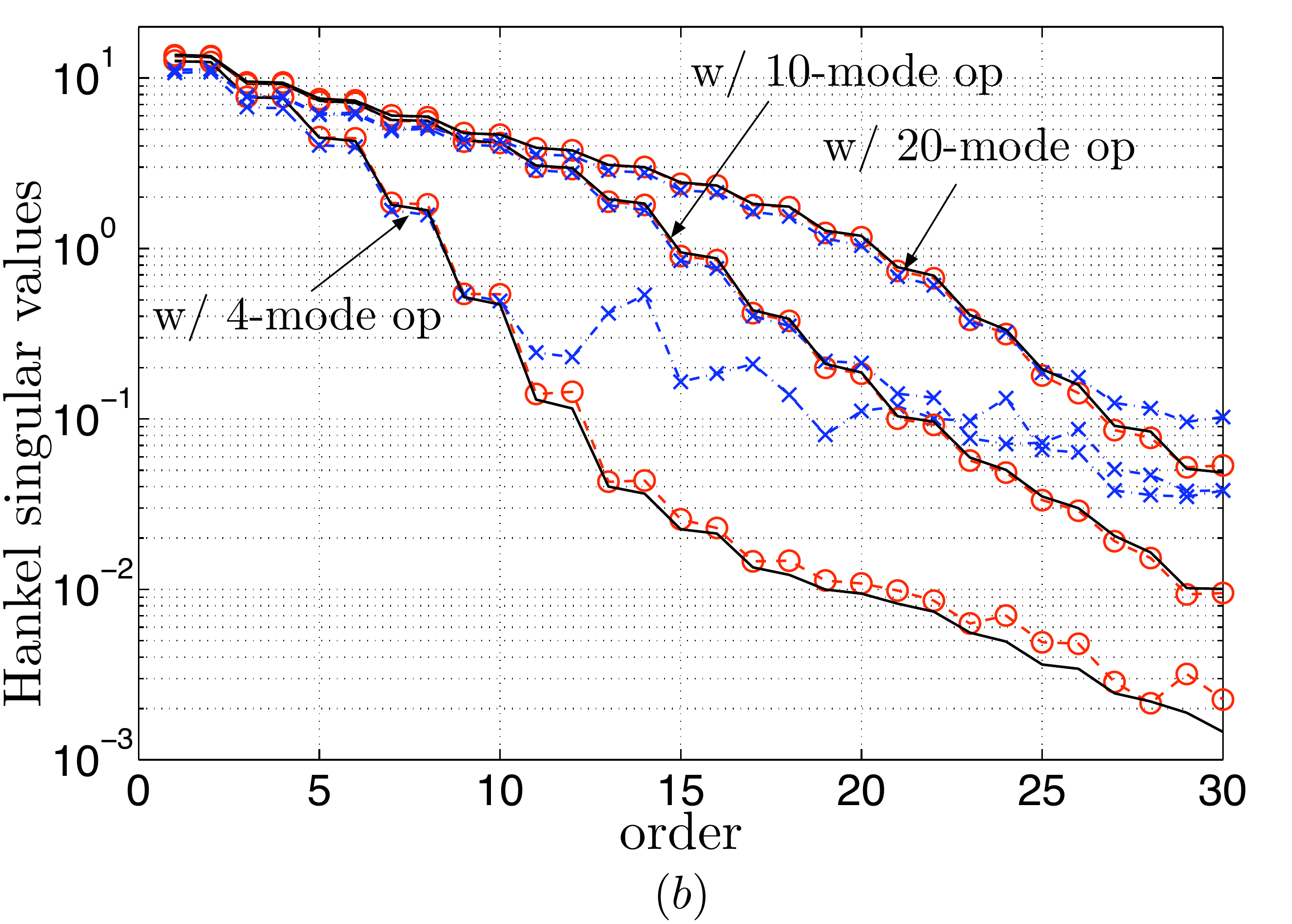}\\
\includegraphics[width=2.5in]{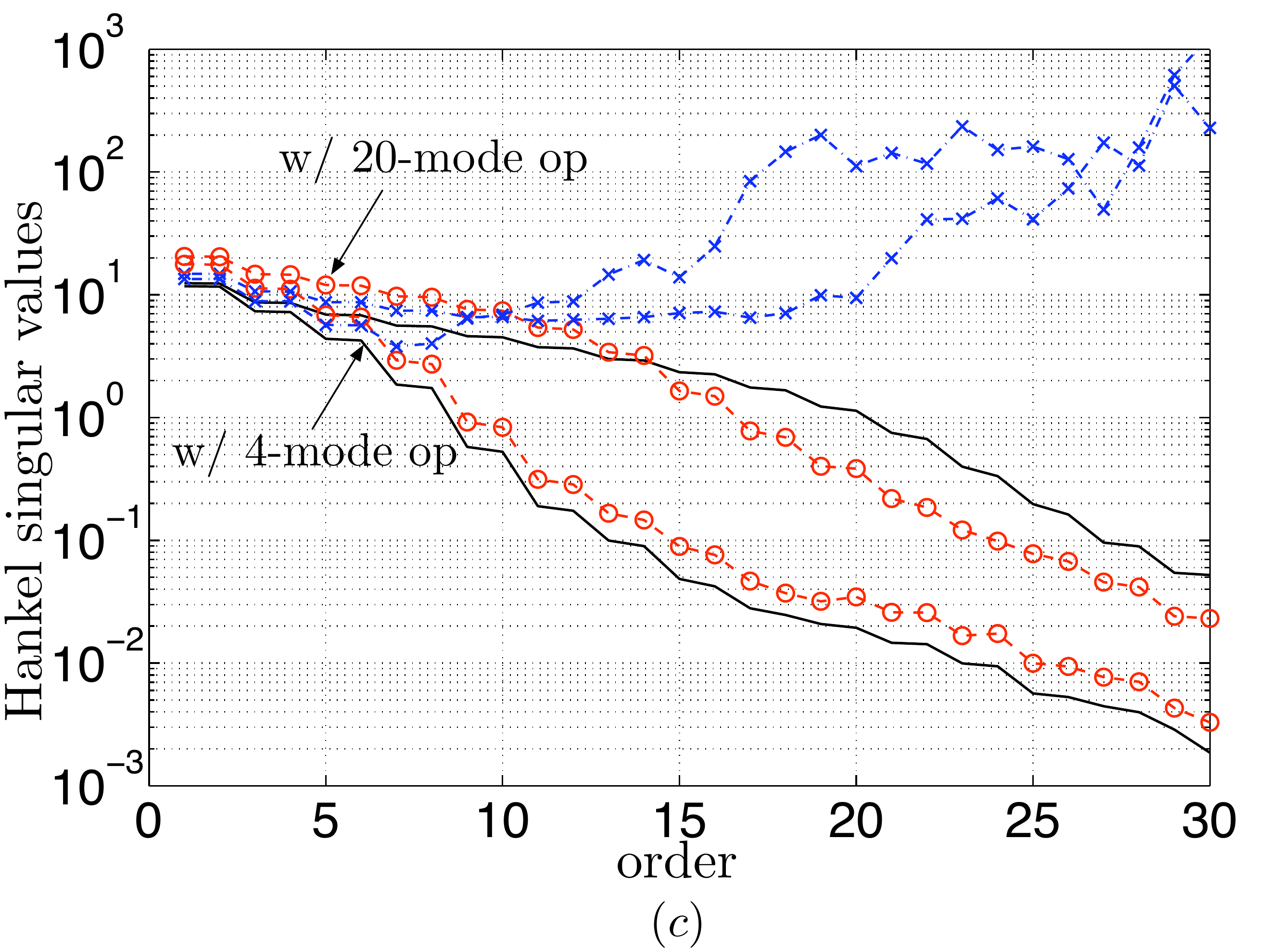}
\caption{Comparison of Gramians computed using (a) ERA, (b) balanced POD, and (c) ERA with pseudo-adjoint modes: The empirical Hankel singular values ( $\solid$) and the diagonal elements of the controllability ({\color{red}$\dashed$, $\circ$}) and observability ({\color{blue}$\chndot$, $\times$}) Gramians with different order of modes (e.g., 4, 10, 20) in output projection.
}
\label{figure:gramians}
\end{figure}

Figure~\ref{figure:gramians} shows the diagonal values of the controllability and observability Gramians, as well as the empirical Hankel singular values, for reduced-order models obtained from three different methods: ERA, balanced POD, and ERA with pseudo-adjoint modes. The models obtained using ERA are more accurate in the sense that the three sets of curves are almost indistinguishable, for all orders of output-projection. However, for balanced POD, the diagonal values of the observability Gramians are accurate only for certain leading modes, the number of which depends on and increases with the order of output projection. This inaccuracy can be attributed to a slight inaccuracy in the adjoint formulation, which in turn results from an approximation in the multi-domain approach used to treat farfield boundary conditions in the immersed boundary method of~\cite{ColTai-07}; see~\cite{Ahuja:2009} for more details. Thus, ERA is advantageous as it does not need any adjoint simulations and results in more balanced Gramians. On the other hand, ERA with pseudo-adjoint modes generates poorly balanced controllability and observability Gramians, as shown in Figure~\ref{figure:gramians}(c).
This is because the leading primal modes and adjoint modes are supported very differently in the spatial domain, and thus the pseudo-adjoint modes, based on linear combination of leading primal modes, poorly approximate the true adjoint modes.   

\subsection{Model performance}
\label{subsec:performance}
\begin{figure}
\centering
\includegraphics[width=0.8\linewidth]{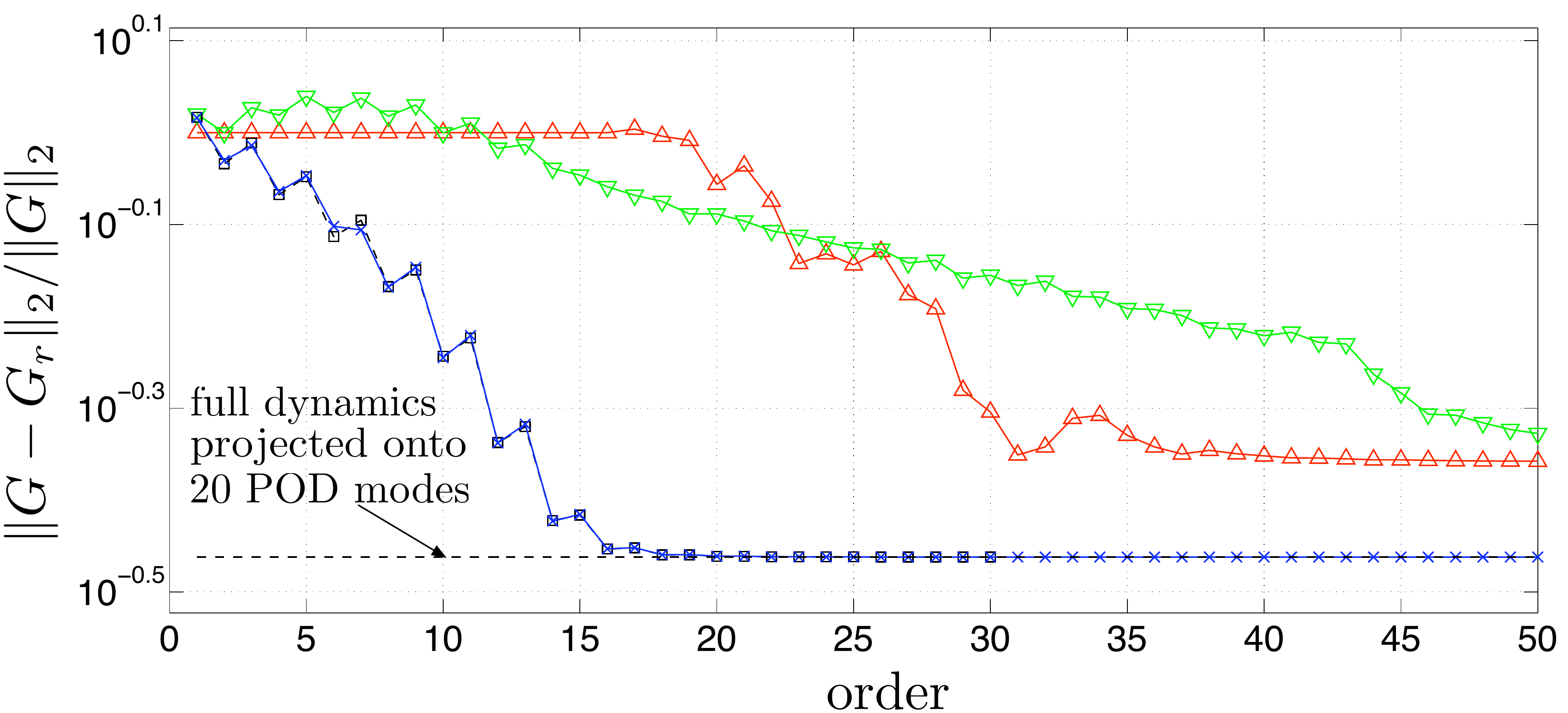}
\caption{$H_2-$norm of the error with increasing order of the reduced-order models: exact output of the output-projected system~($\dashed$); models obtained using balanced POD~({\color{black}$\dashed$,~$\square$}), ERA({\color{blue}$\solid$,~$\times$}), ERA with pseudo-adjoint modes({\color{green}$\solid$,~$\triangledown$}) , and POD({\color{red}$\solid$,$\triangle$}). A 20-mode output projection is used in ERA, balanced POD, and ERA with pseudo-adjoint modes. }
\label{figure:error2norm}
\end{figure}

We can quantify the performance of the various reduced-order models by computing error norms. One such measure is the
$2$-norm of the error between the impulse response of the full linearized system, denoted~$G(t)$, and that of a reduced order model with $r$~modes, denoted by~$G_r(t)$.  We first compute the $2$-norm of the error between the full system (with the entire velocity field as output) and the output-projected system of order~20, shown as the horizontal dashed line in Figure~\ref{figure:error2norm}.  This is the lower error bound for any reduced order model of the given output-projected system.  Results shown in Figure~\ref{figure:error2norm} indicate that the first several low-order models obtained by ERA and balanced POD generate slightly different $2$-norms of error, presumably because of the slight inaccuracy in the adjoint, mentioned previously. For most orders, however, they agree, and both error norms converge to the lower bound as the order of the model increases.  By running more simulation tests, we observe that with higher-order output projections, ERA and balanced POD error norms converge to each other faster when the order of model increases. 

Figure~\ref{figure:error2norm} also shows the 2-norm error plots  for models by ERA with pseudo-adjoint modes, using 20-mode output projection, and  for models computed using standard POD.  Errors of models by ERA with pseudo-adjoint modes converge to the lower bound much slower than ERA/balanced POD.  Errors of models by POD do not start converging until more than nearly 20 modes are used, and they converge to a larger error bound than ERA/balanced POD, again because POD models do not capture the input-output dynamics as well as balanced truncation based models.  

\begin{figure}
\centering
\includegraphics[width=0.8\linewidth]{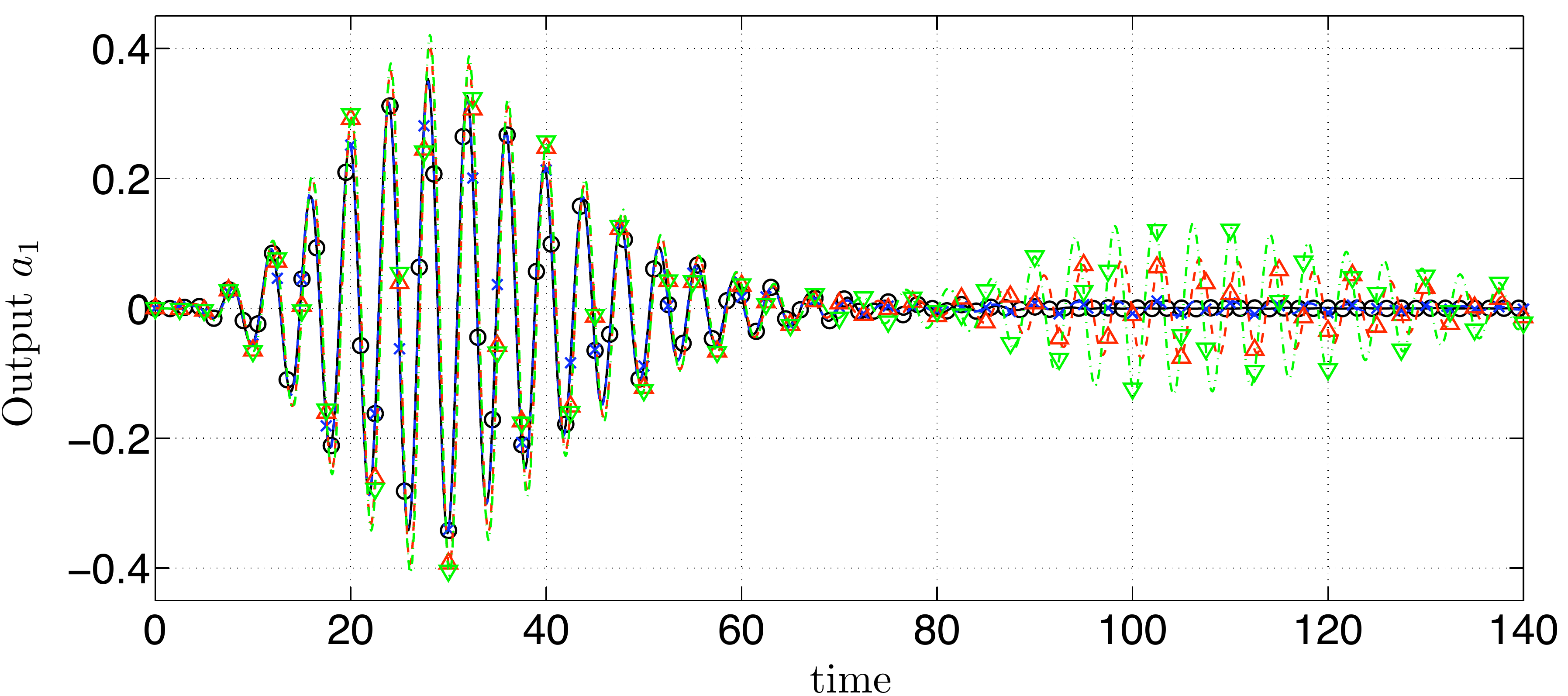}
\caption{The first output, output~$a_1$, from the impulse-response simulation: results of full-simulation
({\color{black}$\solid$,~$\circ$}),  compared with those of $16$-mode reduced order model by ERA({\color{blue}$\dashed$,~$\times$}),  $30$-mode model by ERA with pseudo-adjoint modes({\color{green}$\chndot$ ,~$\triangledown$}) , and $30$-mode model by POD({\color{red}$\dashed$,$\triangle$}) .  A $20$-mode output projection is used in ERA and ERA with pseudo-adjoint modes.}
\label{figure:output} 
\vspace{-0.1in}
\end{figure}

In the time domain, a comparison of the transient response to an impulsive disturbance is shown in Figure~\ref{figure:output}, in which the first output of the reduced-order model is plotted, for a 16-mode model determined by ERA, and for 30-mode models by POD and ERA with pseudo-adjoint modes .  The  16-mode ERA model already accurately predicts the response for all times.   The higher-dimensional, 30-mode models using POD and pseudo-adjoint modes are both stable, and perform reasonably well; however, they over-predict the response, particularly after time $t\approx 80$.  

\begin{figure}
\centering
\includegraphics[width=0.6\linewidth]{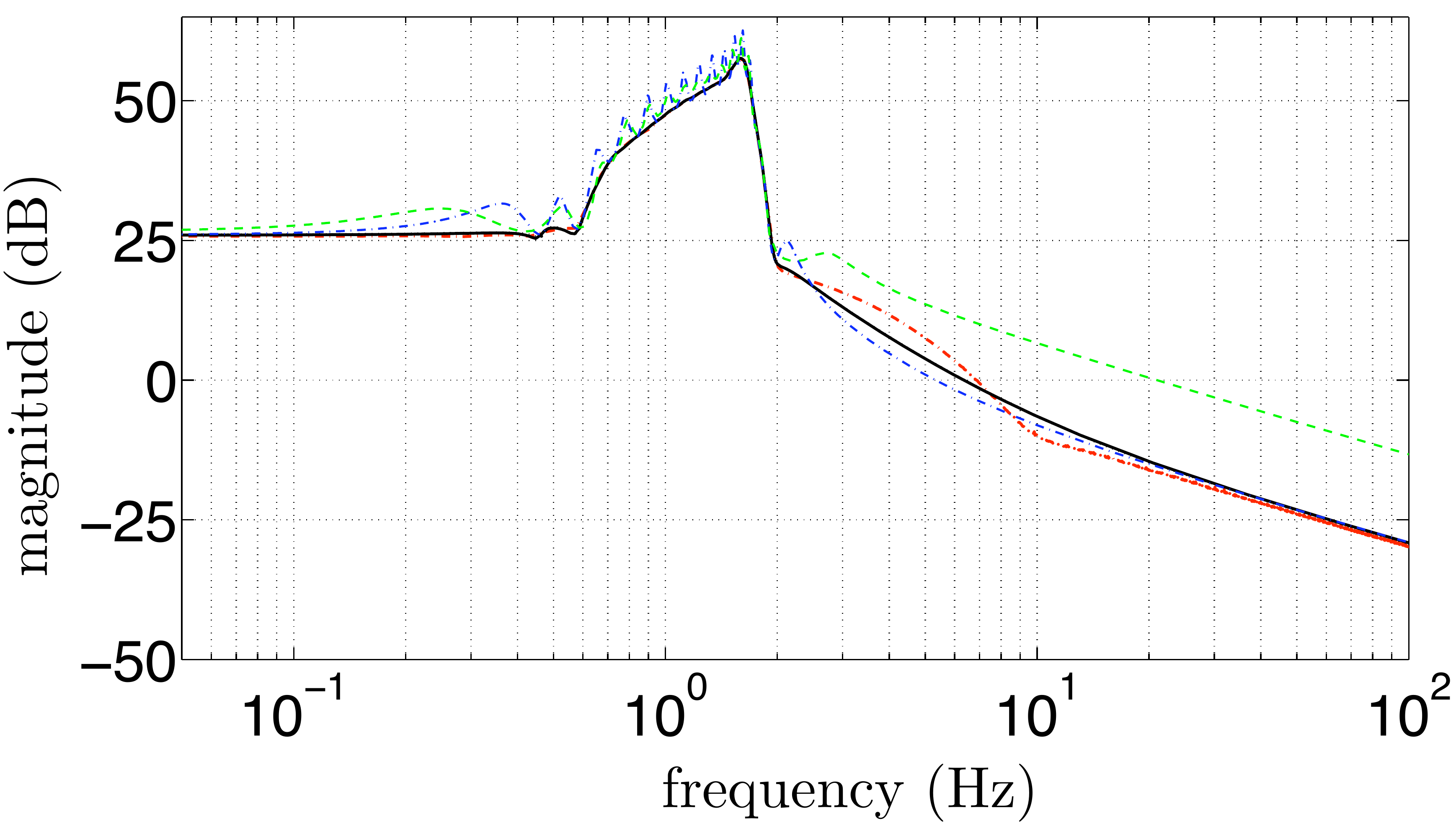} 
\caption{Singular-value plots: 
The full system~({\color{red}$\chndot$}) and 30-mode models obtained using balanced POD~($\dotted$), ERA~($\solid$), ERA with pseudo-adjoint modes({\color{green}$\dashed$}), and POD({\color{blue}$\chndot$}), all with a 20-mode output projection. 
ERA and balanced POD models generate almost identical plots.  
}
\label{figure:sigmaplot}
\end{figure}

We also compare the frequency response of reduced-order models to that of the full system, or more precisely, the full output-projected system.  One way to represent the response of a single-input multiple-output system is by a singular-value plot, a plot of the maximum singular value of the transfer function matrix as a function of frequency.   To generate this plot, a very long simulation of~$5 \times 10^5$~time steps for the full system is performed,  with a random input sampled from a uniform distribution in the range~$(-0.5, 0.5)$.  The output snapshots are projected onto leading POD modes.  The magnitude of the transfer function is then computed from the cross spectrum of the input and outputs (using the Matlab command \texttt{tfestimate}).  Finally, singular-value plots for the full output-projected systems are obtained, with a typical case shown as the red line in Figure~\ref{figure:sigmaplot}.

A typical set of singular-value plots of different reduced order models are presented in Figure~\ref{figure:sigmaplot}.   Results shown in the figure indicate that ERA and balanced POD $30$-mode models, are almost identical, and are close to the corresponding full output-projected system.  
In comparison, Figure~\ref{figure:sigmaplot} also shows sigma plots for $30$-mode models by ERA with pseudo-adjoint modes and by POD.  Note that for computational feasibility,  here the output of the POD model is the first twenty reduced states, i.e., the full-dimensional output of  the POD model are projected onto the leading twenty POD modes.   The frequency responses of the models by POD and ERA with pseudo-adjoint modes capture the resonant peak, but do not match well for frequencies far away from the resonant peak.  These two models both generate spurious peaks in the frequency range of $[0.1, 2]$.

\section{Discussion}
\label{sec:conclusions}
We report that, theoretically, the eigensystem realization algorithm (ERA) and snapshot-based approximate balanced truncation (balanced POD) produce exactly the same reduced order models.   This equivalence implies that ERA balances a pair of  approximate Gramians and thus can be regarded as an approximate balanced truncation method.  Compared to  balanced POD, the main features of ERA are that it does not require data from adjoint systems and therefore can be used with experimental data; furthermore, its construction of the generalized Hankel matrix is computationally an order-of-magnitude cheaper than balanced POD.
Numerical results indicate that ERA can be more accurate than balanced POD in practice, since there can be slight inaccuracies in the adjoint operator used with balanced POD. Balanced POD does have its own advantages, however: unlike ERA, it produces sets of bi-orthogonal modes that are useful for other purposes. Nonlinear models may be obtained by projection onto these modes; and parameters such as Reynolds number can be retained in the reduced-order models generated using these modes.  Balanced POD has also been generalized for unstable systems.  We also examine a modified ERA approach in which one constructs sets of bi-orthogonal modes without using adjoint information, using a matrix pseudo-inverse, as in~\cite{OrSpeCar-08}.  Although this approach provides sets of  bi-orthogonal modes (primal/pseudo-adjoint modes), in general it can not be regarded as an approximate balanced truncation method, since it does not balance the approximate Gramians.

We have demonstrated the methods on an model problem consisting of a disturbance interacting with the flow past an inclined flat plate.  As expected, balanced POD models perform nearly identically to ERA models.  The small differences result because the adjoint simulation required for balanced POD is not a perfect adjoint at the discrete level.  Both procedures work significantly better than standard POD models, or ERA models using pseudo-adjoint modes for projection.

Finally, we emphasize that throughout, we have considered only stable, linear models.  Possible future directions of this work include a generalization to unstable systems, and ultimately to nonlinear systems.

\section{Acknowledgements}
The authors gratefully acknowledge the support by  U. S. Air Force Office of Scientific Research grants FA9550-05-1-0369 and FA9550-07-1-0127.  The authors thank Kunihiko Taira for help with the numerical solver.

\bibliographystyle{abbrv}

\end{document}